\documentclass[12pt]{amsart}

\author{Stefano Decio}
\address{Department of Mathematical Sciences, Norwegian University of Science and Technology, 7491 Trondheim, Norway}
\email{stefano.decio@ntnu.no}

\usepackage{amsmath}
\usepackage{amssymb}
\usepackage{amsfonts}
\usepackage{enumerate}
\usepackage{mathrsfs}
\usepackage{MnSymbol}
\usepackage{mathtools}

\newcommand{\R}{{\mathbf R}}

\newcommand{\ld}{{\lambda}}
\newcommand{\defeq}{\vcentcolon=}

\newtheorem{theorem}{Theorem}
\newtheorem{lemma}{Lemma}
\newtheorem{claim}{Claim}
\newtheorem{corollary}{Corollary}
\newtheorem{proposition}{Proposition}
\newtheorem{thmx}{Theorem}

\theoremstyle{remark}
\newtheorem*{remark}{Remark}
\theoremstyle{definition}

\newtheorem{definition}{Definition}

\DeclareMathOperator{\diver}{div}
\DeclareMathOperator{\Vol}{Vol}
\begin{document}

\begin{abstract}
    We study nodal sets of Steklov eigenfunctions in a bounded domain with $\mathcal{C}^2$ boundary. Our first result is a lower bound for the Hausdorff measure of the nodal set: we show that for $u_{\ld}$ a Steklov eigenfunction, with eigenvalue $\ld\neq 0$, $\mathcal{H}^{d-1}(\{u_{\ld}=0\})\geq c_{\Omega}$, where $c_{\Omega}$ is independent of $\ld$. We also prove an almost sharp upper bound, namely $\mathcal{H}^{d-1}(\{u_{\ld}=0\})\leq C_{\Omega}\ld\log(\ld+e)$.
\end{abstract}

\title{Hausdorff measure bounds for nodal sets of Steklov eigenfunctions}
\title[]{Hausdorff measure bounds for nodal sets of Steklov eigenfunctions}
\maketitle

\section{Introduction}
Let $\Omega$ a bounded domain in $\R^d$, where $d\geq 2$. A Steklov eigenfunction $u_{\ld}\in H^1(\Omega)$ is a solution of
\begin{align}
\label{problem}
    \begin{cases}
    \Delta u_{\ld}=0 \qquad \: \text{in} \ \Omega,\\
    \partial_{\nu}u_{\ld}=\ld u_{\ld} \quad \ \text{on} \ \partial \Omega. \end{cases}
\end{align}
Here and throughout the paper we denote by $\partial_{\nu}$ the outward normal derivative. A number $\ld$ for which a solution to \eqref{problem} exists is called a Steklov eigenvalue, and it is well known that Steklov eigenvalues form a discrete sequence accumulating to infinity. Moreover, Steklov eigenvalues coincide with the eigenvalues of the Dirichlet-to-Neumann operator, which is the operator that maps a function on $\partial\Omega$ to the normal derivative of its harmonic extension in $\Omega$, and a Steklov eigenfunction restricted to $\partial\Omega$ is an eigenfunction of the Dirichlet-to-Neumann operator. For a survey on the Steklov problem outlining many results and open questions see \cite{GP}.\\

Inspired by a famous conjecture of Yau on the Hausdorff measure of nodal sets of Laplace eigenfunctions, an analogous question has been asked for nodal sets of Steklov eigenfunctions (it is stated explicitly in \cite{GP}, for example); the conjecture can be formulated both for interior and boundary nodal sets. For the interior nodal set, the question is as follows:\\
Is it true that there exist positive constants $c,C$, depending only on $\Omega$, such that 
\begin{align}
\label{yauint}
    c\ld\leq \mathcal{H}^{d-1}(\{u_{\ld}=0\}) \leq C\ld ?
\end{align}
Similarly, for the boundary nodal set (which is the nodal set of an eigenfunction of the Dirichlet-to-Neumann operator) one can ask:\\
Is it true that there exist positive constants $c^{'},C^{'}$, depending only on $\Omega$, such that 
\begin{align}
\label{yaubound}
    c'\ld\leq \mathcal{H}^{d-2}(\{u_{\ld}=0\}\cap \partial \Omega) \leq C'\ld ?
\end{align}
Here we do not deal with question \eqref{yaubound} and just note that the upper bound was proved in \cite{Ze} when $\partial\Omega$ is real-analytic. About question \eqref{yauint}, a polynomial upper bound was proved in \cite{GR}, following the corresponding polynomial upper bound in the Laplace-Beltrami eigenfunction case proved in \cite{L2}. On real-analytic surfaces (that is, real-analytic metric in the interior and real-analytic boundary), the full conjecture \eqref{yauint} was established in \cite{PST}. Again in the real-analytic category, the upper bound was recently obtained in any dimension in \cite{Z2}. Concerning lower bounds, as far as we know, the best result was contained in \cite{SWZ}, where the bound $\mathcal{H}^{d-1}(\{u_{\ld}=0\})\geq c\ld^{(2-d)/2}$ is obtained for $\Omega$ a domain with $\mathcal{C}^{\infty}$ boundary (actually, a smooth Riemannian manifold with smooth boundary). The first contribution of the present article is an improvement on the lower bound; we show that the Hausdorff measure of the interior nodal set is bounded below by a constant independent of $\ld$ (so the result is really an improvement over \cite{SWZ} if $d\geq 3$).
\begin{theorem}
\label{lower}
Let $\Omega$ be a bounded domain in $\R^d$ with $\mathcal{C}^2$-smooth boundary, and let $u_{\ld}$ be a solution of \eqref{problem} in $\Omega$, $\ld\neq 0$. Then there exists a constant $c_{\Omega}>0$ independent of $\ld$ such that 
\begin{align}
    \mathcal{H}^{d-1}(\{u_{\ld}=0\})\geq c_{\Omega}.
\end{align}
\end{theorem}
In the previous work \cite{De} we established a density property of the zero set near the boundary, under weaker hypothesis on the boundary regularity: we transcribe the result below.
\begin{thmx}
\label{density}
Let $\Omega$ be a Lipschitz domain in $\R^d$, $d\geq 2$, and let $u_{\ld}$ be a solution of \eqref{problem}, where we assume $\ld\neq 0$. There exists a constant $C=C(\Omega)$ such that 
\begin{align}
\label{den}
    \{u_{\ld}=0\}\cap B\neq \emptyset
\end{align}
for any ball $B$ in $\R^d$ of radius $C/\ld$ centered at a point in $\partial\Omega$. 
\end{thmx}
The proof of Theorem \ref{lower} involves a combination of Theorem \ref{density} and the recent breakthrough by Logunov on Yau's conjecture, \cite{L}. We cannot apply the results of \cite{L} directly and have to do some work to modify the necessary arguments. The fact that we are one power of $\ld$ away from the optimal result is a consequence of the deficiency of the density result, which we can only prove very close to the boundary, and not of the second ingredient. 

\begin{remark}
It will be apparent from the proof that Theorem \ref{lower} extends without much difficulty to the case of manifolds equipped with a $\mathcal{C}^2$-smooth Riemannian metric and $\mathcal{C}^2$ boundary.
\end{remark}

The conjectured upper bound in \eqref{yauint} would be sharp, as the example of a ball shows; the second main contribution of this article is an almost sharp upper bound for Euclidean domains with $\mathcal{C}^2$ boundary. 

\begin{theorem}
\label{upper}
Let $\Omega$ be a bounded domain in $\R^d$ with $\mathcal{C}^2$-smooth boundary, and let $u_{\ld}$ be a solution of \eqref{problem} in $\Omega$. Then there is a constant $C_{\Omega}>0$ independent of $\ld$ such that 
\begin{align}
    \mathcal{H}^{d-1}(\{u_{\ld}=0\})\leq C_{\Omega}\ld\log(\ld+e).
\end{align}
\end{theorem}

\begin{remark}
The proof of Theorem \ref{upper} uses the sharp bounds of Donnelly and Fefferman bound (\cite{DF}) in the interior of the domain and a multiscale induction argument at the boundary, which is based on a version of the Hyperplane Lemma of \cite{L2} and \cite{LMNN}. While, as remarked above, the proof of the lower bound can be extended almost \emph{verbatim} to smooth Riemannian manifolds with boundary, for Theorem \ref{upper} we rely heavily on the fact that $\Omega$ is a Euclidean domain, or at least we have to require that the metric inside $\Omega$ is real analytic; this is because the results of \cite{DF} require real analyticity. Our theorem lies in between previous results on the upper bound: the multiscale argument at the boundary allows for $\mathcal{C}^2$-regularity of the boundary only, as opposed to real analyticity as in the aforementioned paper \cite{Z2}; on the other hand, if the metric inside is assumed to be only $\mathcal{C}^2$ (or $\mathcal{C}^{\infty}$), the best result attainable with these methods is still the polynomial upper bound of \cite{GR}.
\end{remark}

\subsection*{Plan of the paper} We prove Theorem \ref{lower} in Sections \ref{sec:lower} and \ref{sec:nad}; in Section \ref{sec:lower} we discuss a procedure for extending a Steklov eigenfunction across the boundary, which gives rise to an auxiliary equation for which a statement very similar to Logunov's theorem \cite{L} holds (see Theorem \ref{nadirashvili} below), and we use this together with Theorem \ref{density} to prove the lower bound. Section \ref{sec:nad} is quite long and contains the proof of Theorem \ref{nadirashvili}, which requires us to review Logunov's argument carefully and use a combination of classical elliptic estimates and frequency function techniques. Section \ref{sec:upper} is dedicated to the proof of Theorem \ref{upper}.

\section{Lower bound on nodal sets}
\label{sec:lower}
Here we deduce Theorem \ref{lower} using Theorem \ref{density} and ideas stemming from Logunov's solution of a conjecture of Nadirashvili on nodal sets of harmonic functions \cite{L}. In order to do this, we transform a solution to \eqref{problem} into a solution of an elliptic equation in the interior of a domain. To the best of our knowledge, this idea was introduced first in \cite{BL}, and then also applied successfully in \cite{GR}, \cite{Z1}. \\

We now describe this extension procedure, which requires $\partial\Omega$ to be of class $\mathcal{C}^2$; we follow \cite{BL} very closely. There is a $\delta>0$ such that the map $\partial\Omega\times (-\delta,\delta) \ni(y,t)\to y+t\nu(y)$ is one-to-one onto a neighbourhood of $\partial \Omega$ in $\R^d$. We set $d(x)=\text{dist}(x,\partial\Omega)$ and for $\rho\leq \delta$ we define $\Omega_{\rho}=\{x\in \Omega : d(x)<\rho\}$, $\Omega^{'}_{\rho}=\{x\in \R^d : d(x)<\rho\}\setminus \overline{\Omega}$. Let now $u_{\ld}$ be a solution of \eqref{problem}, and for $x\in \Omega_{\delta}\cup \partial\Omega$ define 
\begin{align}
    v(x)=u_{\ld}(x)\exp(\ld d(x));
\end{align}
an easy computation shows that $v$ satisfies:
\begin{align*}
    \begin{cases}
    \diver(A\nabla v)+b(x)\cdot\nabla v+c(x)v=0 \quad \text{in} \ \Omega_{\delta},\\
    \partial_{\nu}v=0 \qquad \qquad \qquad \qquad \quad \qquad \quad \ \text{on} \ \partial \Omega,
    \end{cases}
\end{align*}
where $A=I$, $b=-2\ld\nabla d$, $c=\ld^2-\ld\Delta d$. Consider now the reflection map $\Psi:\Omega_{\delta}\to \Omega^{'}_{\delta}$ given by $\Psi(y+t\nu(y))=y-t\nu(y)$, where $y\in \partial \Omega$; since $v$ satisfies a Neumann boundary condition on $\partial\Omega$, we can extend it "evenly" across the boundary, i.e. set $v(\Psi(x))=v(x)$ for $x\in \Omega_{\delta}$. Denote $\Psi(x)=x'$. Another easy computation shows that on $\Omega^{'}_{\delta}$ the extended function (which we still call $v$) satisfies the equation 
\begin{align*}
    \diver(\widetilde{A}\nabla v)+\widetilde{b}\cdot \nabla v + \widetilde{c}v=0,
\end{align*}
where $\widetilde{A}(x')=\nabla \Psi(x)(\nabla \Psi(x))^T$, $\widetilde{b}^{i}(x')=-\sum_{j}\partial_{x'_j}\widetilde{a}^{ij}(x')+\Delta \Psi^i(x)+\nabla\Psi^i(x)\cdot b(x)$, $\widetilde{c}(x')=c(x)$. Consider now $D=\Omega_{\delta} \cup \partial\Omega \cup \Omega^{'}_{\delta}$; we abuse notation and denote by $A,b,c$ the functions that are equal to the previous $A,b,c$ in $\Omega_{\delta}$ and equal to $\widetilde{A},\widetilde{b},\widetilde{c}$ in $\Omega^{'}_{\delta}$. In \cite{BL} it is shown that $A$ is Lipschitz across $\partial\Omega$ with Lipschitz constant depending only on $\Omega$, and $A$ is uniformly positive definite, again with constant depending only on $\Omega$. Pasting together the pieces, one obtains that $v$ is a strong solution of the uniformly elliptic equation 
\begin{align}
\label{aux}
    \diver(A\nabla v)+b\cdot \nabla v + cv=0
\end{align}
 in $D$, with $A$ Lipschitz, $\|A\|_{L^{\infty}(D)}\leq C$, $\|b\|_{L^{\infty}(D)}\leq C\ld$ and $\|c\|_{L^{\infty}(D)}\leq C\ld^2$. \\
 
 We want to study equation \eqref{aux} at wavelength scale. In order to deal with its zero set we use the theorem below, which is just an extension to more general equations of the aforementioned theorem of Logunov on harmonic functions (\cite{L}); its proof, which merely consists of a tedious but necessary verification that Logunov's argument carries over in this slightly more general setting, is relegated to the next section. We warn the reader that below and in the rest of the paper we do not explicitly indicate dependence of the constants on the dimension. 

\begin{theorem}
\label{nadirashvili}
Consider a strong solution of the equation
\begin{align}
    \label{ell}
    Lu=\diver(A\nabla u)+b\cdot \nabla u + cu=0
\end{align}
in $B=B(0,1)\subset\R^d$, with the following assumptions on the coefficients:
\begin{enumerate}[(i)]
    \item A is a uniformly positive definite matrix, that is $A(x)\xi\cdot\xi\geq \alpha |\xi|^2$ for any $\xi\in \R^d$;
    \item A is Lipschitz, that is $\sum_{i,j}|a^{ij}(x)-a^{ij}(y)|\leq \gamma|x-y|$;
    \item $\sum_{i,j}\|a^{ij}\|_{L^{\infty}(B)}+\sum_{i}\|b^{i}\|_{L^{\infty}(B)}\leq K$;
    \item $c\geq 0$ and $\|c\|_{L^{\infty}(B)}\leq \varepsilon_0$, where $\varepsilon_0$ is a small enough constant depending on $\alpha, \gamma, K$. 
\end{enumerate}
Then there exist $r_0=r_0(\alpha,\gamma,K)<1$, $c_0=c_0(\alpha,\gamma,K)$ such that for any solution $u$ of \eqref{ell} and any ball $B(x,r)\subset B(0,r_0)$ for which $u(x)=0$ we have the lower measure bound:
\begin{align}
\label{nad}
    \mathcal{H}^{d-1}(\{u=0\}\cap B(x,r))\geq c_0r^{d-1}.
\end{align}
\end{theorem}
 
 Assume now that $\ld$ is large enough depending on $\Omega$ and consider a ball $B(x_0,\varepsilon/\ld)\subset D$, where $\varepsilon$ is a small enough constant, with smallness depending only on $\Omega$. We set $v_{x_0,\ld}(x)=v(x_0+\varepsilon x/\ld)$ for $x\in B=B(0,1)$; note that $v_{x_0,\ld}$ satisfies the equation 
 \begin{align}
 \label{loc}
     \diver(A_{x_0,\ld}\nabla v_{x_0,\ld})+b_{x_0,\ld}\cdot \nabla v_{x_0,\ld} + c_{x_0,\ld}v_{x_0,\ld}=0,
 \end{align}
where the ellipticity constant of $A_{x_0,\ld}$ is the same as that of $A$ and the Lipschitz constant is the same if not better, and the coefficients satisfy the bounds $\|A_{x_0,\ld}\|_{L^{\infty}(B)}\leq C$, $\|b_{x_0,\ld}\|_{L^{\infty}(B)}\leq C\varepsilon$ and $\|c_{x_0,\ld}\|_{L^{\infty}(B)}\leq C\varepsilon^2$. Note that if $\ld$ is large enough then $c_{x_0,\ld}\geq 0$. If we then take $\varepsilon$ small enough, $v_{x_0,\ld}$ satisfies equation \eqref{ell} and assumptions \textit{(i)-(iv)} with constants $\alpha,\gamma,K$ depending only on $\Omega$. By Theorem \ref{density}, any ball centered at $\partial\Omega$ of radius $C/\ld$ contains a zero of the Steklov eigenfunction $u_{\ld}$, and hence of $v$. We can reduce the radius of the balls and take a maximal disjoint sub-collection of balls $B(x_i,C_1/\ld)\subset D$, $x_i\in \overline\Omega$, such that $v(x_i)=0$, and consider the corresponding rescaled functions $v_{x_i,\ld}$; we can assume that $C_1<r_0$, so that by Theorem \ref{nadirashvili} we obtain that
\begin{align}
    \mathcal{H}^{d-1}(\{v_{x_i,\ld}=0\}\cap B(0,C_1))\geq cC_1^{d-1}.
\end{align}
Note also that $\mathcal{H}^{d-1}(\{u_{\ld}=0\}\cap B(x_i,C_1/\ld)\cap \Omega)\sim \mathcal{H}^{d-1}(\{v=0\}\cap B(x_i,C_1/\ld))\sim \varepsilon^{d-1}\ld^{1-d}\mathcal{H}^{d-1}(\{v_{x_i,\ld}=0\}\cap B(0,C_1))\geq \widetilde{C}\ld^{1-d}$, where $\widetilde{C}$ depends on $\Omega$ only. Since there are $\sim \ld^{d-1}$ such balls $B(x_i,C_1/\ld)$, we obtain that 
\begin{align*}
    \mathcal{H}^{d-1}(\{u_{\ld}=0\})\geq c_{\Omega}
\end{align*}
and Theorem \ref{lower} is proved. 

\begin{remark}
If one could improve the result of Theorem \ref{density} by showing that every ball of radius $C/\ld$ centered at any point in a corona of fixed (independent of $\ld$) size around the boundary contains a zero of $u_{\ld}$, the optimal lower bound $\mathcal{H}^{d-1}(\{u_{\ld}=0\})\gtrsim \ld$ would follow immediately by the preceding argument (actually more easily, since one could directly apply Logunov's result without the need to go through Theorem \ref{nadirashvili}). 
\end{remark}

\section{Proof of Theorem \ref{nadirashvili}}
\label{sec:nad}
This entire long section is dedicated to the proof of Theorem \ref{nadirashvili}. We follow essentially the arguments of \cite{L}, which carry through in this setting with few changes; the difference is that we have to use more general elliptic estimates, such as a weaker form of maximum principle, and a frequency function that takes into account the lower order terms in the equation. In \ref{subsec:ellest} and \ref{subsec:freq} we introduce the main tools we need in the proof, namely classical elliptic estimates and the monotonicity of the frequency function. Subsection \ref{subsec:outline} will serve as a break from technicalities: here we try to convey an idea of the scheme of the proof to the reader. Subsections \ref{subsec:asymmetry} to \ref{subsec:finalproof} contain the actual body of the proof with full details.\\

Throughout the section we consider the operator $L$ defined by \eqref{ell} satisfying conditions \textit{(i)-(iv)}. It will be convenient to denote by $L_1=L-cI$ the operator without the zeroth order term. 
\subsection{Elliptic estimates} 
\label{subsec:ellest}
We first recall some standard elliptic estimates for $L$, paraphrasing the results in \cite{GT} in our notation. Note that whenever we consider a bounded domain we can assume for our purposes that it is contained in the unit ball, so we can ignore the dependency of the constants on the diameter of $\Omega$, and on the radius of balls contained in $\Omega$. We start with the weak maximum principle. 
\begin{theorem}[\cite{GT}, Theorem 9.1]
\label{weakmax1}
Let $L_1 u\geq -\delta$ in a bounded domain $\Omega$. Then 
$$\sup_{\Omega} u \leq \sup_{\partial\Omega} u^{+}+C|\delta|,$$
where $C=C(\alpha,\gamma,K)$.
\end{theorem}
\begin{corollary}
\label{weakmax2}
Let $Lu=0$ in a bounded domain $\Omega$, with $\varepsilon_0$ in \textit{(iv)} small enough. Then 
\begin{align}
\label{max}
    \sup_{\Omega} u \leq 2\sup_{\partial\Omega} u^{+}
\end{align}
\begin{proof}
We can assume $\sup_{\Omega} u\geq 0$. Since $Lu=0$, we have that $L_1 u=-cu\geq -\varepsilon_0\sup_{\Omega} u$ using assumption \textit{(iv)}. By Theorem \ref{weakmax1} then $\sup_{\Omega} u \leq \sup_{\partial\Omega} u^{+}+C\varepsilon_0\sup_{\Omega} u$, and the corollary follows as soon as $C\varepsilon_0\leq 1/2$.
\end{proof}
\end{corollary}
The next theorem is a local pointwise estimate for subsolutions.
\begin{theorem}[\cite{GT}, Theorem 9.20]
\label{supest}
Let $Lu\geq-\delta$ in $\Omega$. Then for any ball $B(x,2R)\subset \Omega$ and any $p>0$ we have
\begin{align}
\label{sup}
    \sup_{B(x,R)}u\leq C_1\left\{\strokedint_{B(x,2R)}(u^+)^p\right\}^{1/p}+C_2|\delta|,
\end{align}
where $C_1$ and $C_2$ depend on $\alpha,K$ and $p$.
\end{theorem}
\begin{remark}
    In Theorem 9.20 in \cite{GT}, the constants depend on $R$. However they get worse as $R$ increases and improve as $R$ decreases; in this work we will only be concerned with small $R$, so that we can ignore the dependency on it. 
\end{remark}
We now come to the weak Harnack inequality and then the full Harnack inequality.
\begin{theorem}[\cite{GT}, Theorem 9.22]
\label{weakharnack}
Let $Lu\leq \delta$ in $\Omega$, and suppose that $u$ is non-negative in a ball $B(x,2R)\subset \Omega$. Then 
\begin{align}
\label{weakhar}
    \left\{\strokedint_{B(x,R)}u^p\right\}^{1/p}\leq C\left(\inf_{B(x,R)}u+|\delta|\right),
\end{align}
where $p$ and $C$ are positive numbers depending on $\alpha$ and $K$.
\end{theorem}

\begin{theorem}[\cite{GT}, Corollary 9.25]
\label{fullharnack}
Let $Lu=0$ in $\Omega$,and suppose that $u$ is non-negative in a ball $B(x,2R)\subset \Omega$. Then 
\begin{align}
\label{har}
    \sup_{B(x,R)}u\leq C\inf_{B(x,R)}u,
\end{align}
where $C=C(\alpha,K)$.
\end{theorem}

\begin{corollary}
Let $Lu=0$ in $\Omega$. If $u(x_0)\geq 0$ and $B(x_0,R)\subset \Omega$, then the inequality 
\begin{align}
\label{har2}
    \sup_{B(x_0,\frac{2}{3}R)}|u|\leq C\sup_{B(x_0,R)}u
\end{align}
holds for $C=C(\alpha,K)$.
\begin{proof}
Call $M=\sup_{B(x_0,R)}u$ and consider the function $h=M-u$, which is non-negative in $B(x_0,R)$. Note that $Lh=cM$, so that $|Lh|\leq \varepsilon M$. By applying to $h$ in order Theorem \ref{supest} and Theorem \ref{weakharnack} with $\delta=\varepsilon M$, one gets that 
\begin{align*}
    \sup_{B(x_0,\frac{2}{3}R)}(M-u)\leq C_1\left\{\strokedint_{B(x_0,\frac{3}{4}R)}u^p\right\}^{1/p}+C_2\varepsilon M \\ \leq C_3 \inf_{B(x_0,\frac{3}{4}R)}(M-u)+C_4\varepsilon M \leq C_5 M,
\end{align*}
where the last inequality holds because $u(x_0)\geq 0$. Hence we obtain $\sup_{B(x_0,\frac{2}{3}R)}(-u)\leq CM$. Since clearly we have that $\sup_{B(x_0,\frac{2}{3}R)}u\leq M$, the corollary is proved.
\end{proof}
\end{corollary}

\subsection{Frequency function and doubling index} 
\label{subsec:freq}
The frequency function, which as far as we know was used first by Almgren and then subsequently developed in the works of Garofalo and Lin (see \cite{GL1},\cite{GL2}), is a powerful tool in the study of unique continuation and zero sets of elliptic PDEs. We are now going to define it for operators of the form \eqref{ell} and state some of its properties, following mainly \cite{GL2} and \cite{HL}. 

Let $u\in W^{1,2}_{loc}(B)$ be a solution of \eqref{ell}. In \cite{GL2} and \cite{HL} a metric $g(x)=\sum_{i,j}g_{ij}(x)dx_i\otimes dx_j$ is introduced in the following way: let first  
\begin{align*}
    \overline{g_{ij}}(x)=a^{ij}(x)(\det A)^{\frac{1}{d-2}},
\end{align*}
where, as customary, $a^{ij}$ denote the entries of the matrix $A^{-1}$. To define $\overline{g_{ij}}$ we assume here $d\geq 3$; if $d=2$, we can just add a 'mute' variable. Next, one defines
\begin{align*}
    r(x)^2=\sum_{i,j}\overline{g_{ij}}(0)x_ix_j; \quad \quad \eta(x)=\sum_{k,l}\overline{g^{kl}}(x)\frac{\partial r}{\partial {x_k}}(x)\frac{\partial r}{\partial {x_l}}(x).
\end{align*}
Finally, one sets $$g_{ij}(x)=\eta(x)\overline{g_{ij}}(x).$$ Note that $\eta$ is a positive Lipschitz function with Lipschitz constant depending on $\alpha,\gamma, K$. Let $G$ be the matrix $(g_{ij})$ and denote $|g|=\det(G)$. We can now write equation \eqref{ell} as 
\begin{align*}
    \diver_g(\mu(x)\nabla_g u)+b_g(x)\cdot \nabla_g u+c_g(x)u=0,
\end{align*}
where $\mu=\eta^{-\frac{d-2}{2}}$ is a Lipschitz function in $B$ with $C_1\leq \mu(x)\leq C_2$, $b_g=Gb/\sqrt{|g|}$, $c_g=c/\sqrt{|g|}$. Note that, since $|g|^{-1/2}$ is a Lipschitz function bounded above and below by constants depending on $\alpha,\gamma, K$ only, $b_g$ and $c_g$ satisfy analogous bounds to $b$ and $c$ in \eqref{ell}. The following quantities are then introduced, where the integrals are with respect to the measure induced by the metric $g$:
\begin{align*}
    &H(x,r)=\int_{\partial B(x,r)}\mu u^2;
    \\ &D(x,r)=\int_{B(x,r)}\mu |\nabla_g u|^2;
    \\ &I(x,r)=\int_{B(x,r)}\mu |\nabla_g u|^2+ub_g\cdot\nabla_g u + c_gu^2.
\end{align*}
The frequency function is finally defined as 
\begin{align}
\label{freq}
    \beta(x,r)=\frac{2rI(x,r)}{H(x,r)}.
\end{align}
Compared with the definition in \cite{GL2},\cite{HL} there is an extra factor of $2$ for aesthetic reasons in later formulas. More often than not, we will forget about the point $x$ and only write the dependance on the radius $r$. The key property of the frequency function is the following almost monotonicity:
\begin{theorem}
\label{monotonicity}
There are constants $r_0,c_1,c_2$ depending on $\alpha,\gamma,K$ such that 
\begin{align}
\label{mon}
    \beta(x,r)\leq c_1+c_2\beta(x,r_0)
\end{align}
for $r\in (0,r_0)$. Moreover, $c_2$ can be chosen to be $1+\varepsilon$ for any $\varepsilon>0$ if $r_0=r_0(\varepsilon)$ is small enough.
\end{theorem}
\begin{remark}
The statement of Theorem \ref{monotonicity} is implicit in \cite{GL2}, and the proof is contained there; in \cite{HL} the theorem is stated as here and the proof given is essentially the one of \cite{GL2}. The second assertion is not explicitly stated in \cite{GL2} or \cite{HL} and needs some justification. In \cite{GL2} and \cite{HL}, the strategy to prove the theorem is the following: one defines $\Omega_{r_0}=\{r\in (0,r_0) : \beta(r)>\max(1,\beta(r_0))\}$ and proves that it is an open subset of $\R$ and therefore it can be decomposed as $\Omega_{r_0}=\cup_{j=1}^{+\infty}(a_j,b_j)$ with $a_j$ and $b_j$ not belonging to $\Omega_{r_0}$; it is then showed that $\beta'(r)/\beta(r)\geq -C$ for any $r\in\Omega_{r_0}$. By integration, one has that $\beta(r)\leq \beta(b_j)\exp\{C(b_j-r)\}$ for any $r\in (a_j,b_j)$. Since $b_j\notin \Omega_{r_0}$, this implies that the constant $c_2$ can be chosen to be $\exp{(Cr_0)}$, which is close to $1$ if $r_0$ is small.
\end{remark}

In the course of the proof of Theorem \ref{monotonicity} in \cite{GL2} and \cite{HL} the differentiation formula 
\begin{align*}
    H'(r)=\left(\frac{d-1}{r}+O(1)\right)H(r)+2I(r)
\end{align*}
is obtained; the formula can be rewritten as
\begin{align}
\label{diff}
    \frac{d}{dr}\left(\log\frac{H(r)}{r^{d-1}}\right)=O(1)+\frac{\beta(r)}{r}.
\end{align}
The next statement is an immediate consequence of it.
\begin{proposition}
\label{increasing}
    There is a constant $C$ depending on $\alpha,\gamma,K$ such that the function ${e^{Cr}H(r)}/{r^{d-1}}$ is increasing for $r\in(0,r_0)$.
\end{proposition}
From \eqref{diff} and almost monotonicity \eqref{mon}, by integration one obtains the following:
\begin{proposition}
    The two-sided inequality 
    \begin{align}
    \label{scaling}
        c\left(\frac{r_2}{r_1}\right)^{c_2^{-1}\beta(r_1)-c_3}\leq \frac{H(r_2)}{H(r_1)}\leq C 
        \left(\frac{r_2}{r_1}\right)^{c_2\beta(r_2)+c_3}
    \end{align}
holds, where again $c_2$ can be chosen to be $1+\varepsilon$ if $r_0$ is small enough.
\end{proposition}
From now on we denote with letters $c,C,c_1\dots$ constants which may vary from line to line and that depend only on $\alpha,\gamma,K$ without explicitly saying so every time. Additional dependencies will be indicated. We now define a quantity related to the frequency function, the doubling index.
\begin{definition}
For $B(x,2r)\subset B$, the doubling index $\mathcal{N}(x,r)$ is defined by 
\begin{align}
\label{doubdef}
   2^{\mathcal{N}(x,r)}=\frac{\sup_{B(x,2r)}|u|}{\sup_{B(x,r)}|u|}.
\end{align}
\end{definition}
The doubling index and the frequency function are comparable in the following sense:
\begin{lemma}
\label{comparability}
Let $\varepsilon>0$ be sufficiently small, and let $r_0$ be so small that the constant $c_2$ in \eqref{scaling} is $1+\varepsilon$; then, for $4r<r_0$,
\begin{align*}
    \beta(x,r(1+\varepsilon))(1-100\varepsilon)-c\leq \mathcal{N}(x,r)\leq \beta(x,2r(1+\varepsilon))(1+100\varepsilon)+c.
\end{align*}
\end{lemma}
The proof of Lemma \ref{comparability} is an easy computation using the elliptic estimate \eqref{sup}, Proposition \ref{increasing} and \eqref{scaling}; in fact, by \eqref{sup}, 
\begin{align*}
    \sup_{B(x,r)}|u|^2\leq C_{\varepsilon}\strokedint_{B(x,(1+\varepsilon) r)}|u|^2, 
\end{align*}
and further 
\begin{align*}
    \strokedint_{B(x,(1+\varepsilon) r)}|u|^2\leq C H((1+\varepsilon)r)/r^{d-1}
\end{align*}
by integration and Proposition \ref{increasing}. From here on the computation is identical to the one in Lemma 7.1 in \cite{L2}. Using this, one can derive a scaling property for the doubling index; see Lemma 7.2 and Lemma 7.3 in \cite{L2} for details of the computation. 
\begin{lemma}
    Given any $\varepsilon\in (0,1)$, there exist $r_0(\varepsilon)>0$ and $C(\varepsilon)>0$, such that for $u\in W^{1,2}(B)$ a solution of \eqref{ell} and any $0<2r_1\leq r_2\leq r_0$ we have 
    \begin{align}
    \label{doubscal}
        \left(\frac{r_2}{r_1}\right)^{\mathcal{N}(x,r_1)(1-\varepsilon)-C}\leq \frac{\sup_{B(x,r_2)}|u|}{\sup_{B(x,r_1)}|u|}\leq \left(\frac{r_2}{r_1}\right)^{\mathcal{N}(x,r_2)(1+\varepsilon)+C}.
    \end{align}
\end{lemma}
As a consequence, also the doubling index is almost monotonic in the sense that 
\begin{align*}
    \mathcal{N}(x,r_1)(1-\varepsilon)-C\leq \mathcal{N}(x,r_2)(1+\varepsilon)+C.
\end{align*}

\subsection{An informal outline of the proof}
\label{subsec:outline}
We include here a brief discussion of the scheme of the proof avoiding details and technicalities; the latter are all included in the next subsections. Let us first note that in dimension two Theorem \ref{nadirashvili} is an easy consequence of the weak maximum principle (Corollary \ref{weakmax2}): if $u$ vanishes at the center of a ball, the weak maximum principle tells us that there can be no small loops of zeros containing the center and therefore the nodal component containing the center must exit the ball, implying that its length must be greater than the diameter of the ball. \\
In higher dimension, this simple argument does not give any lower measure bound because a priori the nodal set could be a very thin tube crossing the ball. However, a slightly more sophisticated argument, still using essentially only the maximum principle, does give a non-optimal lower bound: we prove in Proposition \ref{naive} below that if $u(x)=0$,
\begin{align*}
    \mathcal{H}^{d-1}(\{u=0\}\cap B(x,r))\geq cr^{d-1}N^{2-d},
\end{align*}
where $N$ is an upper bound for the doubling index $\mathcal{N}(x,r/2)$. Note that when $d=2$ this is already optimal, as it should be. If $d\geq 3$, this naive lower bound gets worse as the doubling index gets larger. This however contradicts intuition, since we are dealing with solutions of elliptic PDEs: if the doubling index is large, meaning that there is strong growth of $u$, then there should be many zeros. This suggests that one could use induction on $N$ to promote the naive lower bound to the optimal one. The key to achieving this is Proposition \ref{key} below, which shows that if the doubling index is comparable to $N>>1$ in balls of radii $r/4$ to $r$ (we call this 'stable growth', see Definition \ref{stablegrowth} below), there are many zeros in the ball of radius $r$; more precisely, there are at least $[\sqrt{N}]^{d-1}f(N)$, with $f(N)\to\infty$ as $N\to\infty$, disjoint balls of radius $r/\sqrt{N}$ such that $u$ vanishes at the center. The fact that $f(N)$ grows with $N$ essentially shows that indeed there are more zeros as the doubling index increases, and it is needed to close the induction in \ref{subsec:finalproof}. \\
The proof of Proposition \ref{key} uses crucially Theorem \ref{count2}, which tells us that if a cube is partitioned into some large number $B^d$ of subcubes, the number of subcubes which have doubling indices dropping by an amount increasing with $B$ compared to the doubling index of the original cube form the vast majority of the subcubes. The argument goes as follows: since the doubling index is comparable to $N$ on scales $r/4$ to $r$, we can assume that in the ball of radius $r/4$, $|u|\leq 1$, while in the ball of radius $r/2$, $|u|\geq 2^{cN}$. We then connect points where $u$ is small to points where $u$ is large by many chains of cubes (called 'tunnels' later): since there is considerable growth of $u$ from one endpoint of the tunnel to the other, the Harnack inequality tells us that there must be zeros and the growth happens in the cubes with zeros; an application of Theorem \ref{count2} gives us that most of the cubes in the tunnel have doubling index much smaller than $N$, so that the growth from one endpoint to the other cannot be realized in very few cubes, and hence each tunnel must have many cubes with zeros. The formal proof is a matter of quantifying what 'small', 'large', 'few' and 'many' mean. \\
The only issue remanining is ensuring that there are balls of stable growth: this is done in Claim \ref{smaller}, and the proof uses the estimates in \ref{subsec:shell} which are consequences of the almost monotonicity of the frequency function. \\

Let us emphasize once again that the proof scheme described above is due to Aleksandr Logunov, and it appeared first in \cite{L}. In our case we have to adapt it to elliptic equations with lower order terms, but the more general estimates that we need are collected above in subsections \ref{subsec:ellest} and \ref{subsec:freq} and using those estimates the proof runs in the same way as for harmonic functions. 

\subsection{Local asymmetry} 
\label{subsec:asymmetry}
We now derive a lower estimate for the relative volume of the set $\{u>0\}$ in balls centered at zeros of $u$, and consequently a non-optimal lower estimate for the measure of the zero set. The estimate and the proof are analogous to the Laplace-Beltrami eigenfunctions case, for which see for example \cite{LM} and \cite{M}. For the reader's convenience, we reproduce here essentially the same proof as \cite{LM}. 

\begin{proposition}
\label{naive}
Let $B(x,r)\subset B$, and $u$ be a solution of \eqref{ell} such that $u(x)=0$. Suppose that $\mathcal{N}(x,r/2)\leq N$, where $N$ is a positive integer. Then the lower measure bound 
\begin{align}
\label{easy}
    \mathcal{H}^{d-1}(\{u=0\}\cap B(x,r))\geq cr^{d-1}N^{2-d}
\end{align}
holds for some $c>0$.
\begin{proof}
For notational simplicity we assume $x=0$ and denote $B_r=B(0,r)$. We can also safely assume that $N\geq 4$, say. Note that by \eqref{har2} and \eqref{max} we have that $\sup_{B_{r/2}}|u|\leq C \max_{\partial B_{3r/4}}u$, so that 
\begin{align*}
    \frac{\max_{\partial B_{r}}u}{\max_{\partial B_{3r/4}}u}\leq C_1 \frac{\sup_{B_{r}}|u|}{\sup_{B_{r/2}}|u|}\leq C_1 2^N.
\end{align*}
Let now $r_j=r(3/4+j/4N)$ for $j=0,1,\dots,N$, and consider the concentric spheres $S_j=\{|x|=r_j\}$. Denote by $m_j^+=\max_{S_j}u$ and $m_j^-=\min_{S_j}u$. From the weak maximum principle \eqref{max} (applied to $u$ as well as $-u$) it follows that $m_j^+>0$, $m_j^-<0$, $m_j^+\leq 2 m_{j+1}^+$ and $|m_j^-|\leq 2 |m_{j+1}^-|$. For $j=0,1,\dots, N-1$, denote by $\tau_j^+=m_{j+1}^+/m_j^+$, $\tau_j^-=|m_{j+1}^-|/|m_j^-|$; from the above, $\tau_j^{+/-}\geq 1/2$. Moreover we have that 
\begin{align*}
    \tau_0^+\dots\tau_{N-1}^+=\frac{\max_{\partial B_{r}}u}{\max_{\partial B_{3r/4}}u}\leq C_1 2^N,
\end{align*}
so at most $N/4$, say, of the $\tau_j^+$ are greater than some $C$ independent of $N$. The same holds for the $\tau_j^-$, so that for at least $N/2$ indices there holds $m_{k+1}^+\leq C m_k^+$ and $|m_{k+1}^-|\leq C |m_k^-|$. Consider each such $k$ and let $x_0\in S_k$ be such that $u(x_0)=m_k^+$. Denote by $b$ the ball centered at $x_0$ of radius $r/8N$; then by \eqref{max} and the choice of $k$
\begin{align*}
    \sup_b u \leq \sup_{\{|x|\leq r_{k+1}\}} u \leq 2m_{k+1}^+\leq C m_k^+.
\end{align*}
Applying \eqref{har2}, we then get that $\sup_{b/2}|u|\leq C m_k^+$. We now use this last inequality and the elliptic gradient estimate (see for instance \cite{GT}, Theorem 8.32)
\begin{align*}
    \sup_{B(y,s/2)}|\nabla u|\leq (C/s)\sup_{B(y,s)}|u|
\end{align*} 
for $y=x_0$ and $s=r/16N$ to get, for $x\in B(x_0,\theta r/N)$ where $\theta$ is a sufficiently small number, 
\begin{align*}
    u(x)\geq u(x_0)-|x-x_0|\sup_{b/4}|\nabla u|\geq m_k^+-C\theta m_k^+\geq 0.
\end{align*}
We thus found a ball centered on $S_k$ of radius $\theta r/N$ where $u$ is positive, call it $b_+$. Replace now $u$ with $-u$, which is also a solution of \eqref{ell}: repeating the argument above with $m_k^-$ and $\tau_k^-$ instead of $m_k^+$ and $\tau_k^+$ gives us a ball centered on $S_k$ of radius $\theta r/N$ where $u$ is negative, call it $b_-$. Now consider the sections of the two balls with hyperplanes through the origin that contain the center of the balls: any path within the annulus $\{r_{k-1}<|x|<r_{k+1}\}$ that connects these two sections contains a zero of $u$, since $u$ is positive on $b_+$ and negative on $b_-$. This implies that the measure of the zero set is greater than the measure of the section of the balls, that is to say:
\begin{align*}
    \mathcal{H}^{d-1}(\{x: r_{k-1}<|x|<r_{k+1}, u(x)=0\})\geq c\left(\frac{r}{N}\right)^{d-1}.
\end{align*}
The above holds for all indices $k$ for which $m_{k+1}^+\leq C m_k^+$ and $|m_{k+1}^-|\leq C |m_k^-|$, and recall that there are at least $N/2$ such indices. Summing the inequality above over those indices, we see that \eqref{easy} holds. 
\end{proof}
\end{proposition}

\begin{remark}
Note that the argument above also shows that 
\begin{align*}
    \frac{\Vol(\{u>0\}\cap B(x,r))}{\Vol(B(x,r))}\geq \frac{c}{N^{d-1}}
\end{align*}
if $u(x)=0$, which is analogous to the best known lower bound (when $d\geq 3$) for the local asymmetry of Laplace eigenfunctions (\cite{M}).
\end{remark}

\subsection{Counting doubling indices}
\label{subsec:count}
We now recall some very useful results from \cite{L2}, \cite{L}, \cite{LM} that allow to find many small cubes with better doubling index than the original ball (or cube). The proofs are combinatorial in nature. First we define a version of the doubling index for cubes, which are more suitable for partitioning than balls. Given a cube $Q$ and a solution $u$ of \eqref{ell}, we define the doubling index $N(Q)$
as
\begin{align*}
    N(Q)=\sup_{\{x\in Q, r<diam(Q)\}}\log\frac{\sup_{B(x,10dr)}|u|}{\sup_{B(x,r)}|u|}.
\end{align*}
The constant $10d$ is there for technical reasons and the reader should not worry about it. It is clear that with this definition $N(Q_1)\leq N(Q_2)$ if $Q_1\subset Q_2$. Theorem \ref{count} below was proved in \cite{L2}, and then extended in \cite{GR} to the more general equation \eqref{ell}; the proof combines an accumulation of growth result (the Simplex Lemma, Lemma 2.1 in \cite{L2} and Proposition 3.1 in \cite{GR}) and a propagation of smallness result (The Hyperplane Lemma, Lemma 4.1 in \cite{L2} and Proposition 3.2 in \cite{GR}). The Hyperplane Lemma is a consequence of quantitative Cauchy uniqueness, which we state in a simple version below; it can be obtained from a very general result in \cite{ARRV} (Theorem 1.7). See also \cite{Li}.

\begin{proposition}
\label{cauchy}
Let $D$ be a bounded domain with $\mathcal{C}^2$ boundary, and let $B$ be a ball of radius $\rho<1$. Let $u$ be a solution of \eqref{ell} in $D\cap B$, $u\in\mathcal{C}^1(\overline{D}\cap B)$. There exist $\beta=\beta(\alpha,\gamma,K,D,\rho)\in (0,1)$ and $C=C(\alpha,\gamma,K,D,\rho)>0$ such that if $|u|\leq 1$, $|\nabla u|\leq \rho^{-1}$ in $D\cap B$ and $|u|\leq \eta$, $|\nabla u|\leq \eta\rho^{-1}$ on $\partial D\cap B$, where $\eta$ is a real number, then 
\begin{align*}
    |u(x)|\leq C\eta^{\beta}
\end{align*}
for any $x\in D\cap \frac{1}{2}B$.
\end{proposition}

\begin{remark}
In \cite{L2} and \cite{GR} Proposition \ref{cauchy} is applied when $\partial D$ is flat; this is sufficient to prove the theorem below. We will use the proposition in the non-flat case later in Section 4, to prove a version of the Hyperplane Lemma. 
\end{remark}

\begin{theorem}[\cite{L2}, Theorem 5.1 and \cite{GR}, Theorem 4.1]
\label{count}
There exist a constant $c>0$ and an integer $A>1$ depending on the dimension only, and positive numbers $N_0=N_0(\alpha,\gamma,K)$, $R_0=R_0(\alpha,\gamma,K)$ such that for any cube $Q\subset B(0,R_0)$ the following holds:\\
If $Q$ is partitioned into $A^d$ equal subcubes, then the number of subcubes with doubling index greater than $\max\left(N(Q)/(1+c),N_0\right)$ is less than $\frac{1}{2}A^{d-1}$.
\end{theorem}

Starting from Theorem \ref{count}, in \cite{L} an iterated version is proved, which is the one decisively used in the proof of the lower bound on zero sets. We state it below and refer to \cite{L} for the proof. 
\begin{theorem}[\cite{L}, Theorem 5.3]
\label{count2}
There exist positive constants $c_1,c_2,C$ and an integer $B_0>1$ depending on the dimension only, and positive numbers $N_0=N_0(\alpha,\gamma,K)$, $R_0=R_0(\alpha,\gamma,K)$ such that for any cube $Q\subset B(0,R_0)$ the following holds:\\
If $Q$ is partitioned into $B^d$ equal subcubes, where $B>B_0$ is an integer, then the number of subcubes with doubling index greater than $\max\left(N(Q)2^{-c_1\log B/\log\log B},N_0\right)$ is less than $CB^{d-1-c_2}$.
\end{theorem}

\subsection{Estimates in a spherical shell} 
\label{subsec:shell}
In the following we always indicate by $u$ a solution of \eqref{ell}; the frequency function and doubling index are relative to $u$. Consider a ball $B(p,s)\subset B(0,r_0/4)$; we are going to establish some estimates for the growth of $u$ near a point of maximum. Let $x\in\partial B(p,s)$ be a point where 
the maximum of $|u|$ on $\overline{B(p,s)}$ is almost attained, in the sense that $\sup_{B(p,s)}|u|\leq 2|u(x)|$; the existence of such an $x$ is guaranteed by Corollary 1. Call $M=|u(x)|$. In the next two lemmas we will assume that there is a large enough number $N$ and 
$$\delta\in \left(\frac{1}{\log^{100}N},\frac{1}{8}\right)$$
such that 
\begin{align}
\label{stablefreq}
    N/10\leq \beta(p,t) \leq 10^4 N 
\end{align}
for $t\in I\defeq (s(1-\delta),s(1+\delta))$.
\begin{lemma}[variation on Lemma 4.1, \cite{L}]
\label{shell}
    Let \eqref{stablefreq} be satisfied. There exist positive constants $C,c$ such that 
    \begin{align}
        \label{small}
        \sup_{B(p,s(1-\delta))}|u|\leq CM2^{-c\delta N},\\
        \label{large}
        \sup_{B(p,s(1+\delta))}|u|\leq CM2^{C\delta N}.
    \end{align}
    \begin{proof}
    Let us prove \eqref{small} only. By \eqref{scaling} and \eqref{stablefreq}, we have that 
    \begin{align}
    \label{scaling2}
        \left(\frac{t_2}{t_1}\right)^{N/30}\leq \frac{H(p,t_2)}{H(p,t_1)}\leq C 
        \left(\frac{t_2}{t_1}\right)^{10^5N},
    \end{align}
    for $t_1<t_2\in I$, where we assume that $r_0$ is small enough to take $c_2=2$ in \eqref{scaling}. We estimate:
    \begin{align*}
        M^2\geq C_1 s^{-d+1} H(p,s)\geq C_1 s^{-d+1} H(p,s(1-\delta/2))(1+\delta/2)^{N/30},
    \end{align*}
    where the first inequality is just the estimate of the $L^2$-norm by the $L^{\infty}$-norm and the second inequality comes from \eqref{scaling2}. By integration and Proposition \ref{increasing} we have that 
    \begin{align*}
        sH(p,s(1-\delta/2))=s\int_{\partial B(p,s(1-\delta/2))}|u|^2\geq C_2 \int_{B(p,s(1-\delta/2))}|u|^2.
    \end{align*}
    Let now $\widetilde{x}$ be a point on $\partial B(p,s(1-\delta))$ where the sup of $|u|$ on $B(p,s(1-\delta))$ is almost attained, i.e. $\sup_{B(p,s(1-\delta))}|u|\leq 2|u(\widetilde{x})|$, and call $\widetilde{M}=|u(\widetilde{x})|$. Note now that 
    \begin{align*}
        \int_{B(p,s(1-\delta/2))}|u|^2\geq \int_{B(\widetilde{x},\delta s/2)}|u|^2\geq C_3 (\delta s)^d \strokedint_{B(\widetilde{x},\delta s/2)}|u|^2;
    \end{align*}
    moreover, by \eqref{sup} we have that 
    \begin{align*}
        \widetilde{M}^2\leq C_4\strokedint_{B(\widetilde{x},\delta s/2)}|u|^2.
    \end{align*}
    Combining the estimates we obtain 
    \begin{align*}
        M^2\geq C_5 \delta^d(1+\delta/2)^{N/30}\widetilde{M}^2.
    \end{align*}
    Since $\log(1+\delta/2)\geq \delta/4$, it follows easily from the above and $\delta\gtrsim 1/\log^{100}N$ that $M^2\geq C_6\exp(N\delta/100)\widetilde{M}^2$, from which one obtains \eqref{small} recalling the definition of $M$ and $\widetilde{M}$.
    \end{proof}
\end{lemma}

Using the properties of the doubling index, we now derive some estimates on small balls close to $x$; we keep on denoting by $x$ the point on $\partial B(p,s)$ where the maximum of $|u|$ on $\overline{B(p,s)}$ is almost attained. 
\begin{lemma}[variation on Lemma 4.2, \cite{L}]
\label{fine}
    Let \eqref{stablefreq} be satisfied. There exists $C>0$ such that 
    \begin{align}
    \label{d1}
        \sup_{B(x,\delta s)}|u|\leq M 2^{C\delta N+C}
    \end{align}
    and for any $\widetilde{x}$ with $d(x,\widetilde{x})\leq \delta s/4$
    \begin{align}
    \label{d2}
        \mathcal{N}(\widetilde{x},\delta s/4)\leq C\delta N+C,\\
        \label{d3}
        \sup_{B(\widetilde{x},\delta s/10N)}|u|\geq M 2^{-C\delta N\log N-C}.
    \end{align}
    \begin{proof}
    Note that since $B(x,\delta s)\subset B(p,s(1+\delta))$, the first estimate \eqref{d1} is an immediate consequence of \eqref{large}. By definition of doubling index and \eqref{d1} we have that 
    \begin{align*}
        2^{\mathcal{N}(\widetilde{x},\delta s/4)}\leq \frac{\sup_{B(\widetilde{x},\delta s/2)}|u|}{\sup_{B(\widetilde{x},\delta s/4)}|u|}\leq \frac{\sup_{B(x,\delta s)}|u|}{M}\leq 2^{C\delta N+C},
    \end{align*}
    and \eqref{d2} is proved. Now recall the scaling properties \eqref{doubscal}; by those and \eqref{d2} we obtain 
    \begin{align*}
        \frac{\sup_{B(\widetilde{x},\delta s/4)}|u|}{\sup_{B(\widetilde{x},\delta s/10N)}|u|}\leq (40N)^{2\mathcal{N}(\widetilde{x},\delta s/4)+C_1}\leq (40N)^{C_1\delta N+C_1}\\
        \leq 2^{C_2\delta N\log N+C_2\log N}\leq 2^{C_3\delta N\log N+C_3},
    \end{align*}
    where the last inequality holds because $\delta\gtrsim 1/\log^{100}N$. Since, by the distance condition, $\sup_{B(\widetilde{x},\delta s/4)}|u|\geq |u(x)|=M$, \eqref{d3} follows. 
    \end{proof}
\end{lemma}

\subsection{Finding many balls around the zero set} 
\label{sec:manyzeros}
We follow the arguments in Section 6 of \cite{L}, in the reformulation contained in \cite{LM2}; the estimates in the spherical shell will be used together with the combinatorial results on doubling indices. We use the notion of "stable growth", which is taken from \cite{LM2} and was not present in \cite{L}. 
\begin{definition}
\label{stablegrowth}
We say that $u$ has a stable growth of order $N$ in a ball $B(y,s)$ if $\mathcal{N}(y,s/4)\geq N$ and $\mathcal{N}(y,s)\leq 1000N$.
\end{definition}
The number $1000$ does not have any special meaning, it is just a large enough numerical constant. The following result is the key to the proof of the lower bound. 
\begin{proposition}[variation on Proposition 6.1, \cite{L}]
\label{key}
Let $B(p,2r)\subset B(0,r_0)$. There exist a number $N_0>0$ large enough such that for $N>N_0$ and any solution $u$ of \eqref{ell} that has a stable growth of order $N$ in $B(p,r)$, the following holds:\\
There exist at least $[\sqrt{N}]^{d-1}2^{c_1\log N/\log\log N}$ disjoint balls $B(x_i,r/\sqrt{N})\subset B(p,r)$ such that $u(x_i)=0$.
\begin{proof}
Assume without loss of generality that $\sup_{B(p,r/4)}|u|=1$. The stable growth assumption then implies that 
$$\sup_{B(p,r/2)}|u|\geq 2^N \text{ and } \sup_{B(p,2r)}|u|\leq 2^{CN}.$$
We denote by $x$ the point on $\partial B(p,r/2)$ where the maximum over $\overline{B(p,r/2)}$ is almost attained, so that by the above $|u(x)|\geq 2^{N-1}$. We now divide the ball $B(p,2r)$ into cubes $q_i$ of side length $cr/\sqrt{N}$, and organize these cubes into tunnels in the following way: the centers of the cubes in each tunnel lie on a line parallel to the segment that connects $p$ and $x$. A tunnel contains at most $C\sqrt{N}$ cubes. Let us call a cube $q_i$ good if 
\begin{align}
    \label{good}
        N(q_{i})\leq \max\left(\frac{N}{2^{c\log N/\log\log N}},N_0\right)
\end{align}
for some constant $c$. We will call a tunnel good if it contains only good cubes; by Theorem \ref{count2}, most of the cubes are good and most of the tunnels are good. Another application of Theorem \ref{count2} gives the following:
\begin{claim}
\label{manytunnels}
 The number of good tunnels containing at least one cube with distance from $x$ less than $r/\log^2N$ is greater than $c(\sqrt{N}/\log^2 N)^{d-1}$.
\end{claim}
The proof of the proposition is then completed with the help of the next claim.
\begin{claim}
\label{manycubes}
 Any good tunnel that contains at least one cube with distance from $x$ less than $r/\log^2N$ also contains at least $2^{c_2\log N/\log\log N}$ cubes with zeros of $u$.
 \begin{proof}
 Take one such tunnel $T$. Note that $T$ contains at least one cube $q_a\subset B(p,r/4)$, so that $\sup_{q_a}|u|\leq 1$. Call $q_b$ a cube in $T$ with distance from $x$ less than $r/\log^2N$; we want to show that the supremum of $|u|$ over $q_b$ is large. To this end, we apply Lemma \ref{fine} with $\delta\sim 1/\log^2N$ and $\widetilde{x}$ being the center $x_b$ of the cube $q_b$. By the stable growth assumption and the comparability of doubling index and frequency (Lemma \ref{comparability}), \eqref{stablefreq} is satisfied for $N$ large enough. Then \eqref{d3} gives us that 
 \begin{align*}
     \sup_{B(x_b,\delta r/10N)}|u|\geq |u(x)|2^{-CN/\log N-C}
 \end{align*}
 and hence, recalling that $|u(x)|\geq 2^{N-1}$, 
 \begin{align*}
     \sup_{\frac{1}{2}q_b}|u|\geq 2^{cN}.
 \end{align*}
 We now follow $T$ from $q_a$ to $q_b$ and find many zeros. The proof is at this point identical to the one given in \cite{L}; for completeness we provide the details. We enumerate the cubes $q_i$ from $q_a$ to $q_b$ so that $q_a$ is the first and $q_b$ is the last. Since $T$ is a good tunnel, by \eqref{good} we have that for any two adjacent cubes 
 \begin{align*}
     \log\frac{\sup_{\frac{1}{2}q_{i+1}}|u|}{\sup_{\frac{1}{2}q_{i}}|u|}\leq \log\frac{\sup_{4q_{i}}|u|}{\sup_{\frac{1}{2}q_{i}}|u|}\leq \frac{N}{2^{c_3\log N/\log\log N}}.
 \end{align*}
 We split the set of indices $S$ into two sets $S_1$ and $S_2$, where $S_1$ is the set of $i$ such that $u$ does not change sign in $\overline{q_{i}}\cup\overline{q_{i+1}}$ and $S_2=S\setminus S_1$. The advantage of this is the possibility to use the Harnack inequality on $S_1$; the aim is to get a lower bound on the cardinality of $S_2$. In fact, for $i\in S_1$, by \eqref{har} we have that 
 \begin{align*}
     \log\frac{\sup_{\frac{1}{2}q_{i+1}}|u|}{\sup_{\frac{1}{2}q_{i}}|u|}\leq C_1.
 \end{align*}
 We then estimate 
 \begin{align*}
     \log\frac{\sup_{\frac{1}{2}q_{b}}|u|}{\sup_{\frac{1}{2}q_{a}}|u|}=\sum_{S_1}\log\frac{\sup_{\frac{1}{2}q_{i+1}}|u|}{\sup_{\frac{1}{2}q_{i}}|u|}+\sum_{S_2}\log\frac{\sup_{\frac{1}{2}q_{i+1}}|u|}{\sup_{\frac{1}{2}q_{i}}|u|}\\
     \leq |S_1|C_1+|S_2|\frac{N}{2^{c_3\log N/\log\log N}};
 \end{align*}
 on the other hand, recall that 
 \begin{align*}
     \log\frac{\sup_{\frac{1}{2}q_{b}}|u|}{\sup_{\frac{1}{2}q_{a}}|u|}\geq {cN}.
 \end{align*}
 Combining the two estimates one obtains 
 \begin{align*}
     |S_1|C_1+|S_2|\frac{N}{2^{c_3\log N/\log\log N}}\geq {cN}
 \end{align*}
 and noting that $|S_1|C_1\leq C_1\sqrt{N}\leq cN/2$ we conclude 
 \begin{align*}
     |S_2|\geq c_32^{c_3\log N/\log\log N}.
 \end{align*}
 The last quantity is larger than $2^{c_2\log N/\log\log N}$ if $N$ is large enough, and the claim is proved.
 \end{proof}
\end{claim}
 It is now a straightforward matter to finish the proof of Proposition \ref{key}: by Claim \ref{manytunnels} there are at least $c(\sqrt{N}/\log^2 N)^{d-1}$ tunnels satisfying the hypothesis of Claim \ref{manycubes}, hence at least $c(\sqrt{N}/\log^2 N)^{d-1}2^{c_2\log N/\log\log N}$ cubes that contain zeros of $u$; the last quantity can be made larger than $(\sqrt{N})^{d-1}2^{c_1\log N/\log\log N}$, and then one replaces cubes by balls. 
\end{proof}
\end{proposition}

\subsection{Proof of the lower bound}
\label{subsec:finalproof}
We take $r_0$ so small that \eqref{mon}, \eqref{scaling}, Lemma \ref{comparability} and \eqref{doubscal} hold. Denote by $N(0,r_0)=\sup_{\{B(x,r)\subset B(0,r_0)\}}\mathcal{N}(x,r)$; we define 
\begin{align*}
    F(N)=\inf\frac{\mathcal{H}^{d-1}(\{u=0\}\cap B(x,\rho))}{\rho^{d-1}},
\end{align*}
where the $\inf$ is taken over all balls $B(x,\rho)\subset B(0,r_0)$ and all solutions $u$ of \eqref{ell} such that $u(x)=0$ and $N(0,r_0)\leq N$. Theorem \ref{nadirashvili} then follows immediately from the following:
\begin{theorem}
$F(N)\geq c$, where $c$ is independent of $N$.
\begin{proof}
Let $u$ be a solution of \eqref{ell} in competition for the $\inf$ in the definition of $F(N)$; let $F(N)$ be almost attained on $u$, in the sense that 
\begin{align}
\label{inf}
    \frac{\mathcal{H}^{d-1}(\{u=0\}\cap B(x,\rho))}{r^{d-1}}\leq 2F(N)
\end{align}
for some $B(x,r)\subset B(0,r_0)$ with $u(x)=0$. Recall the easy bound \eqref{easy}: 
\begin{align}
\label{easy2}
    \frac{\mathcal{H}^{d-1}(\{u=0\}\cap B(x,r))}{r^{d-1}}\geq \frac{c_1}{\mathcal{N}(x,r/4)^{d-2}}\geq \frac{c_1}{N^{d-2}}.
\end{align}
Estimate \eqref{easy2} already finishes the proof if $\mathcal{N}(x,r/4)$ is bounded uniformly in $N$; let us then argue by contradiction and assume that $\mathcal{N}(x,r/4)$ is large enough. Denote $\widetilde{N}=\mathcal{N}(x,r/4)$ and suppose first that $u$ has stable growth of order $\widetilde{N}$. We can then apply Proposition \ref{key} and find at least $[\sqrt{\widetilde{N}}]^{d-1}2^{c\log \widetilde{N}/\log\log \widetilde{N}}$ disjoint balls $B(x_i,r/\sqrt{\widetilde{N}})\subset B(x,r)$ with $u(x_i)=0$. By definition of $F(N)$, there holds:
\begin{align*}
    \mathcal{H}^{d-1}\left(\{u=0\}\cap B\left(x_i,\frac{r}{\sqrt{\widetilde{N}}}\right)\right)\geq F(N) \left(\frac{r}{\sqrt{\widetilde{N}}}\right)^{d-1}.
\end{align*}
Summing over the inequality over all the balls, we obtain 
\begin{align*}
    \mathcal{H}^{d-1}(\{u=0\}\cap B(x,\rho))\geq \left[\sqrt{\widetilde{N}}\right]^{d-1}2^{c\log \widetilde{N}/\log\log \widetilde{N}}F(N) \left(\frac{r}{\sqrt{\widetilde{N}}}\right)^{d-1};
\end{align*}
the quantity on the right can be made larger than $2F(N)r^{d-1}$ if $\widetilde{N}$ is large enough, which is a contradiction with \eqref{inf}. Therefore we would be done if we knew a priori that $u$ has stable growth of order $\widetilde{N}$ in $B(x,r)$, but this is not necessarily the case; fortunately we can find a smaller ball where $u$ has stable growth. 
\begin{claim}
\label{smaller}
If $\mathcal{N}(x,r/4)$ is large enough, there is a number $N_1\gtrsim\mathcal{N}(x,r/4)$ and a ball $B_1\subset B(x,r)$ with radius $r_1\sim r/\log^2 N_1$ such that $u$ has stable growth of order $N_1/\log^2 N_1$ in $B_1$.
\begin{proof} Let us define a modified frequency function as $$\widetilde{\beta}(p,r)=\sup_{t\in (0,r]}\beta(p,t)+c_1,$$ so that $\widetilde{\beta}(p,r)$ is a positive monotonic increasing function. Note that by \eqref{mon} we have that 
\begin{align*}
    \beta(p,r)\leq \widetilde{\beta}(p,r)\leq c_3+2\beta(p,r),
\end{align*}
 and the rightmost term is less than $3\beta(p,r)$ if $\beta(p,r)\geq c_3$. We use the following:
\begin{claim}[\cite{L}, Lemma 3.1]
\label{monotonic}
    Let $f$ be a non-negative, monotonic non-decreasing function in $[a,b]$, and assume $f\geq e$. Then there exist $x\in[a, (a+b)/2)$ and a number $N_1\geq e$ such that 
    \begin{align*}
        N_1\leq f(t) \leq eN_1
    \end{align*}
    for any $t\in (x-\frac{b-a}{20\log^2f(x)},x+\frac{b-a}{20\log^2f(x)})\subset [a,b].$
\end{claim}
We apply Claim \ref{monotonic} to $\widetilde{\beta}(p,\cdot)$ and hence identify a spherical shell of width $\sim r/\log^2 N_1$ about $s\in (2r/3,3r/4)$ where $\widetilde{\beta}(p,\cdot)$ is comparable to $N_1$. Since $\mathcal{N}(x,r/4)$ is large, by Lemma \ref{comparability} and almost monotonicity $\beta(x,t)$ is large for $t>r/2$ and then also $\beta(x,\cdot)$ is comparable to $N_1$ in the spherical shell. In other words, \eqref{stablefreq} holds with $N_1$ and $\delta\sim 1/\log^2 N$. Let now $y\in\partial B_s$ be a point where the maximum is almost attained, as in Lemmas \ref{shell} and \ref{fine}. Take a ball $B_1$ of radius $\sim s/\log^2 N_1$ such that $\frac{1}{4}B_1\subset B_{s(1-\delta)}$ and $y\in \frac{1}{2}B_1$; then \eqref{small} implies that $$\mathcal{N}(\frac{1}{4}B_1)\geq c\frac{N_1}{\log^2 N_1}$$ and \eqref{large} implies that 
$$\mathcal{N}(B_1)\leq C\frac{N_1}{\log^2 N_1},$$
which means that $u$ has stable growth of order $N_1/\log^2 N_1$ in $B_1$, and the claim is proved. 
\end{proof}
\end{claim}
Claim \ref{smaller} gives an order of stable growth that is again large enough to get a contradiction with \eqref{inf} if $\mathcal{N}(x,r/4)$ and hence $N_1$ is large enough. This means that $\mathcal{N}(x,r/4)$ is bounded from above by some $N_0$ independently of $N$, and therefore by \eqref{inf} and \eqref{easy2} we obtain 
\begin{align}
    F(N)\geq \frac{\mathcal{H}^{d-1}(\{u=0\}\cap B(x,r))}{2r^{d-1}}\geq \frac{c_3}{(N_0)^{d-2}}\geq c,
\end{align}
which concludes the proof of the theorem.
\end{proof}
\end{theorem}

\section{Upper bound}
\label{sec:upper}
Here we give the proof of Theorem \ref{upper}. Throughout this section $\partial\Omega$ is assumed to be of class $\mathcal{C}^2$. As remarked in the introduction, the proof uses the Donnelly-Fefferman bound (\cite{DF}) in the interior of the domain and a multiscale induction argument at the boundary. As will be apparent from the proof, the result with a $\mathcal{C}^{\infty}$-metric inside $\Omega$ would follow from an upper bound for zero sets of elliptic PDEs with smooth coefficients that is linear in the frequency; the best we have thus far is polynomial, \cite{L2}.\\

We introduce now a version of the doubling index that takes into account the boundary. Namely, for $x\in\overline{\Omega}$ and $u\in\mathcal{C}(\overline{\Omega})$ an harmonic function, we let 
\begin{align}
    2^{\mathcal{N}_u^*(x,r)}=\frac{\sup_{B(x,2r)\cap \Omega}|u|}{\sup_{B(x,r)\cap \Omega}|u|}.
\end{align}
Note that if $v$ is the extension across the boundary of the Steklov eigenfunction $u_{\ld}$ as in Section 5 and $dist(x,\partial\Omega)\lesssim 1/\ld$, $r\lesssim 1/\ld$, we have that $\mathcal{N}_{u_{\ld}}^*(x,r)\sim\mathcal{N}_v(x,r)$, where $\mathcal{N}_v(x,r)$ is defined as in \eqref{doubdef}; this will allow us to use the almost monotonicity property \eqref{doubscal}. It was proved in \cite{Z1} (using the extension $v$) that for any $r<r_0(\Omega)$ 
\begin{align}
\label{doubest}
    \mathcal{N}_{u_{\ld}}^*(x,r)\leq C\ld,
\end{align}
mirroring a corresponding statement for Laplace eigenfunctions proved by Donnelly and Fefferman. It will once again be convenient to define a maximal version of the doubling index for cubes; for $Q\subset \R^d$ a cube such that $Q\cap\Omega\neq \emptyset$, we set 
\begin{align*}
    N_u^*(Q)=\sup_{x\in Q\cap\overline{\Omega}, r\leq diam(Q)}\mathcal{N}_u^*(x,r).
\end{align*}

\begin{definition}
We call a Whitney cube in $\Omega$ any cube $Q$ such that $c_1 \text{dist}(Q,\partial\Omega) \leq s(Q) \leq c_2 \text{dist}(Q,\partial\Omega)$, where $s(Q)$ is the side length of $Q$ and $c_1$ and $c_2$ are positive dimensional constants. 
\end{definition}

With this notation, we state the following important result of \cite{DF}.
\begin{theorem}
\label{df}
Let $u$ be a harmonic function in $\Omega$. Then there is $C>0$, independent of $u$, such that 
\begin{align}
\label{sharp}
    \mathcal{H}^{d-1}(\mathcal{Z}_u\cap Q)\leq C(N_u^*(Q)+1)s(Q)^{d-1}
\end{align}
for any Whitney cube $Q$. 
\end{theorem}

From now on, we will denote by $u$ a Steklov eigenfunction with eigenvalue $\ld$. We will first use the theorem above to bound the measure of the zero set of $u$ in the interior, up to a distance from the boundary comparable to $1/\lambda$. We will assume $\ld>\ld_0$. As in the previous section, denote $d(x)=dist(x,\partial\Omega)$; Let $c_0$ be a small constant depending only on $\Omega$. We decompose 
\begin{align*}
    \Omega=In\cup{Mid}\cup{Bd},
\end{align*}
where $In=\{x\in\Omega : d(x)\geq c_0\}$, $Mid=\{x\in\Omega : {c_0}/{\ld}<d(x)< c_0\}$, $Bd=\{x\in\Omega : d(x)\leq {c_0}/{\ld}\}$. It follows easily from Theorem \ref{df} and \eqref{doubest} that
\begin{align}
\label{in}
    \mathcal{H}^{d-1}(\mathcal{Z}_u\cap In)\leq C\ld,
\end{align}
with $C$ depending on $\Omega$ only. The next lemma estimates the contribution of the nodal set in $Mid$.
\begin{lemma}
    There is $C>0$ depending only on $\Omega$ such that 
    \begin{align}
    \label{mid}
        \mathcal{H}^{d-1}(\mathcal{Z}_u\cap Mid)\leq C\ld\log\ld.
    \end{align}
    \begin{proof}
    We set $M_k=\{x\in\Omega : c_02^{k-1}/\ld<d(x)<c_02^{k}/\ld\}$, and we have 
    \begin{align*}
        Mid=\bigcup_{k=1}^{c\log\ld}M_k.
    \end{align*}
    We perform a decomposition of $\Omega$ into Whitney cubes with disjoint interior (the statement that this is possible is usually called the Whitney Covering Lemma). Define
    \begin{align*}
        \mathcal{Q}_k=\{\text{Whitney cubes intersecting } M_k\}.
    \end{align*}
    In the following lines we will denote by $|\cdot|$ both the cardinality of a discrete collection and the Lebesgue measure of cubes; it should cause no confusion. Note that if $Q\in \mathcal{Q}_k$, then 
    \begin{align*}
        |Q| \sim \frac{2^{kd}}{\ld^d};
    \end{align*}
    it follows that $|\mathcal{Q}_k|\lesssim 2^{-kd}\ld^{d-1}$. We can then estimate, using Theorem \ref{df} and \eqref{doubest},
    \begin{align*}
        \mathcal{H}^{d-1}(\mathcal{Z}_u\cap Mid)=\sum_{k=1}^{c\log\ld}\mathcal{H}^{d-1}(\mathcal{Z}_u\cap M_k)\leq \sum_{k=1}^{c\log\ld}\sum_{Q\in \mathcal{Q}_k}\mathcal{H}^{d-1}(\mathcal{Z}_u\cap Q)\\
        \lesssim\ld \sum_{k=1}^{c\log\ld}\sum_{Q\in \mathcal{Q}_k}s(Q)^{d-1}\lesssim \ld\sum_{k=1}^{c\log\ld}|\mathcal{Q}_k|\frac{2^{kd}}{\ld^{d-1}}\lesssim \ld\log\ld,
    \end{align*}
    and the lemma is proved.
    \end{proof}
\end{lemma}

To prove Theorem \ref{upper} the only thing left is to estimate $\mathcal{H}^{d-1}(\mathcal{Z}_u\cap Bd)$. We cover $Bd$ with $\sim\ld^{d-1}$ cubes $q_{\ld}$ centered at $\partial\Omega$ of side length $s(q_{\ld})=4c_0/\ld$; then Theorem \ref{upper} follows from \eqref{doubest} and the following:
\begin{proposition}
\label{linear}
Let $q_{\ld}$ be one of the cubes above, and suppose $N_u^*(4q_{\ld})\leq N$. Then 
\begin{align}
\label{bd}
    \mathcal{H}^{d-1}(\mathcal{Z}_u\cap q_{\ld})\leq C(\Omega)(N+1)s(q_{\ld})^{d-1}.
\end{align}
\end{proposition}

\begin{remark}
In the following we will rescale 
\begin{align}
\label{scale}
    h(x)=u(x/\ld)
\end{align}
so that $q_{\ld}$ becomes a cube $Q$ of side length $s<1$, where $s$ is small enough depending on $\Omega$ but independent of $\ld$, and $h$ satisfies $\Delta h=0$ in $10Q\cap\Omega$, $\partial_{\nu}h=h$ on $\partial\Omega\cap \overline{10Q}$. Note that the doubling index is unchanged under this rescaling. Proposition \ref{linear} will follow from 
\begin{align}
\label{rescaled}
    \mathcal{H}^{d-1}(\mathcal{Z}_h\cap Q)\leq C(\Omega)(N+1).
\end{align}
\end{remark}

The main ingredient in the proof of Proposition \ref{linear} is a version of the Hyperplane Lemma of \cite{L2} with cubes touching the boundary, the proof of which uses quantitative Cauchy uniqueness as stated in Proposition \ref{cauchy}. The proof is very similar to the one contained in \cite{LMNN}, we reproduce it here for the reader's convenience. 

\begin{lemma}
\label{hyperplane}
    Let $h$ be as in \eqref{scale}, and $Q$ as in the remark above a cube of side length $s$. There exist $k, N_0$ large enough depending on $s$ and $\Omega$ such that if $Q\cap\partial\Omega$ is covered by $2^{k(d-1)}$ cubes $q_j$ with disjoint interior centered at $\partial\Omega$ of side length $2^{-k}s$, and $N_h^*(Q)=N>N_0$, then there exist $q_{j_0}$ such that $N_h^*(q_{j_0})\leq N/2$.
    \begin{proof}
    We note first that since $\partial\Omega$ is of class $\mathcal{C}^2$, $h$ is harmonic in $10Q\cap\Omega$ and $\partial_{\nu}h=h$ on $\partial\Omega\cap \overline{10Q}$, we can use the extension-across-the-boundary trick described in Section 4, namely consider $v(x)=e^{d(x)}h(x)$; recall that the coefficients of the second order term in the equation satisfied by $v$ are at least Lipschitz. This gives us access to elliptic estimates that hold up to the boundary for $h$. In particular we will use the gradient estimate:
    \begin{align}
    \label{grad}
        \sup_{B(y,r)\cap\overline{\Omega}}|\nabla h|\lesssim \frac{1}{r} \sup_{B(y,2r)\cap\overline{\Omega}}|h|,
    \end{align}
    where the implied constant depends on $s$ and $\Omega$. Denote now by $x_Q\in\partial\Omega$ the center of the cube $Q$. Consider a ball $B$ centered at $x_Q$ such that $2Q\subset B$, and let $M=\sup_{B\cap\Omega}|h|$. By contradiction, suppose that $N_h^*(q_{j})> N/2$ for any $j$; by definition, this implies that for any $j$ there is $x_j\in q_j\cap\Omega$ and $r_j\leq 2^{-k}\sqrt{d}s=\vcentcolon r_0$ such that $\mathcal{N}_h^*(x_j,r_j)>N/2$. Assuming $N$ large enough, we use \eqref{doubscal} to get 
    \begin{align*}
        \sup_{B(x_j,2r_0)\cap\Omega}|h|\leq (C2^{-k})^{N/10}\sup_{B\cap\Omega}|h|\leq Me^{-cNk}
    \end{align*}
    if $k$ is large enough. Using \eqref{grad}, we get 
    \begin{align*}
        \sup_{B(x_j,r_0)\cap\Omega}|\nabla h|\lesssim \frac{1}{r_0}Me^{-cNk},
    \end{align*}
    with implied constant depending on $s$ and $\Omega$. Note that since $q_j\subset B(x_j,r_0)$ the two estimates above give bounds for the Cauchy data of $h$ on $\partial\Omega\cap Q$. On the other hand if $B'$ is the ball centered at $x_Q$ such that $4B'\subset Q$ we have that $\sup_{2B'\cap\Omega}|h|\leq M$, $\sup_{2B'\cap\Omega}|\nabla h|\lesssim \frac{1}{s}M$. Recalling that $r_0=2^{-k}\sqrt{d}s$, we can then apply Proposition \ref{cauchy} with $\eta=2^{k}e^{-cNk}$ to get 
    \begin{align*}
      \sup_{B'\cap\Omega}|h|\leq C(s,\Omega)2^{\beta k c_d}e^{-c\beta Nk}M.
    \end{align*}
    But then 
    \begin{align*}
        \mathcal{N}_h^*(x_Q,\sqrt{d}s)\geq C_d\log\frac{\sup_{B\cap\Omega}|h|}{\sup_{B'\cap\Omega}|h|}\geq C_d(c\beta Nk-c_d\beta k-C),
    \end{align*}
    and the rightmost term is larger than $N$ if $k$ and $N$ are large enough depending on $s$ and $\Omega$; this is a contradiction with $N_h^*(Q)=N$.
    \end{proof}
\end{lemma}

We are now ready to prove Proposition \ref{linear}, or actually \eqref{rescaled}. The argument is an iteration at the boundary; it originates in \cite{LMNN}. 

\begin{proof}[Proof (of \eqref{rescaled})]
First, we consider again $v(x)=e^{d(x)}h(x)$ and its even extension across the boundary (which we still call $v$). Recall from Section 2 that $v$ satisfies an elliptic PDE with Lipschitz second order coefficients and bounded lower order coefficients. The results of \cite{HS} then apply to this situation. Let $Q$ be any cube with $s(Q)<s_0$ small enough. By Theorem 1.7 of \cite{HS}, we have that 
\begin{align*}
    \mathcal{H}^{d-1}(\mathcal{Z}_v\cap B(x, \rho))\leq CN_v(Q)\rho^{d-1}
\end{align*}
for any ball $B(x,\rho)\subset Q$ where $v(x)=0$ and $\rho<\rho_0(N_v(Q))$. Covering $\mathcal{Z}_h\cap Q$ with balls of such small radius and summing the estimate above over all those balls, it follows that there is a function $\widetilde{A}:\R_+\to\R_+$ such that 
\begin{align}
\label{finite}
    \mathcal{H}^{d-1}(\mathcal{Z}_h\cap Q)\leq \widetilde{A}(N_h^*(Q))s(Q)^{d-1}.
\end{align}
Let now $Q$ be as above a cube centered at $\partial\Omega$ of side $s$, with $s$ small enough depending on $\Omega$. Fix a large number $N_0$; if $N_h^*(Q)<N_0$, \eqref{finite} already implies the result. Otherwise, cover $Q\cap\Omega$ with smaller cubes of side length $2^{-k}s$, where $k=k(\Omega)$ is given by Lemma 8, in the following way: first $Q\cap\partial\Omega$ is covered by cubes $q\in\mathcal{B}$ centered at $\partial\Omega$ with disjoint interior, and then the rest of $Q\cap \Omega$ is covered by cubes $q\in \mathcal{I}$ with $dist(q,\partial\Omega)>cs(q)$ for some constant $c>0$ independent of $k$. Cubes in $\mathcal{B}$ will be called boundary cubes and cubes in $\mathcal{I}$ will be called inner cubes; inner cubes are allowed to overlap, while boundary cubes are not. Denote $N_h^*(Q)=N$. By \eqref{sharp} and almost monotonicity there holds
\begin{align*}
    \mathcal{H}^{d-1}(\mathcal{Z}_h\cap(\cup_{q\in\mathcal{I}}q))\leq C(k)Ns^{d-1}.
\end{align*}
By Lemma \ref{hyperplane}, there is a boundary cube, call it $q_0$, such that $N_h^*(q_0)<N/2$. The other cubes in $\mathcal{B}$ will be enumerated from $1$ to $2^{k(d-1)}-1$. We have that 
\begin{align*}
    \frac{\mathcal{H}^{d-1}(\mathcal{Z}_h\cap Q)}{s^{d-1}}\leq CN+\frac{\mathcal{H}^{d-1}(\mathcal{Z}_h\cap q_0)}{s^{d-1}}+\sum_{j=1}^{2^{k(d-1)}-1}\frac{\mathcal{H}^{d-1}(\mathcal{Z}_h\cap q_j)}{s^{d-1}}.
\end{align*}
We define now 
\begin{align*}
    A(N)=\sup\frac{\mathcal{H}^{d-1}(\mathcal{Z}_h\cap q)}{s(q)^{d-1}},
\end{align*}
where the sup is taken over all harmonic functions $h$ in $2Q$ with $\partial_{\nu}h=h$ on $\partial\Omega\cap 2Q$, $N_h^*(Q)\leq N$ and all cubes $q\subset Q$. By \eqref{finite}, $A(N)<+\infty$. From the inequality above, we get 
\begin{align*}
    A(N)\leq C(k)N+A(N/2)2^{-k(d-1)}+(2^{k(d-1)}-1)A(N)2^{-k(d-1)},
\end{align*}
from which 
\begin{align*}
    A(N)<C(k)N+A(N/2).
\end{align*}
(Beware that $C(k)$ changes value from line to line, and depends also on $\Omega$). Iterating the last inequality until $N/2<N_0$, we obtain 
\begin{align*}
    A(N)<C(k)N+A(N_0)<C(k)(N+1),
\end{align*}
which concludes the proof. 
\end{proof}

Theorem \ref{upper} now follows by combining \eqref{in}, \eqref{mid}, \eqref{bd} and \eqref{doubest}. We believe that the extra $\log\ld$ factor is not necessary and is an artificial feature of the proof; it appears in the proof of \eqref{mid} and it is due to the necessity of getting to cubes of side length $\sim \ld^{-1}$.

\section*{Acknowledgements}
This work owes a lot to the patient guidance of Eugenia Malinnikova; discussing the article with her has contributed greatly to both the contents and the presentation, and her reading of several drafts helped spot mistakes and inaccuracies. Many thanks are due to Aleksandr Logunov for his encouragement, as well as for reading this work at various stages and asking useful questions. I am also very grateful for the encouragement of Iosif Polterovich, whose comments helped improve the presentation of this article.

The work for the present article was started while I was a Visiting Student Researcher at the Department of Mathematics at Stanford University; it is a pleasure to thank the department for the hospitality and the nice working conditions.

The author is supported by Project 275113 of the Research Council of Norway.

\end{document}